\documentclass[10pt]{article}
\usepackage{amssymb}
\usepackage{amsmath}
\usepackage[dvips]{graphics}
\usepackage{epsfig}

\begin{document}
\date{}
\author{\textbf{Vassilis G. Papanicolaou}
\\\\
Department of Mathematics
\\
National Technical University of Athens,
\\
Zografou Campus, 157 80, Athens, GREECE
\\\\
{\tt papanico@math.ntua.gr}}
\title{The Weierstrass $\wp$-function of the hexagonal lattice}
\maketitle
\begin{abstract}
We present some properties of the Weierstrass $\wp$-function associated to the hexagonal (or triangular) lattice. In particular, with the help of
an old theorem of I.N. Baker \cite{B} on the characterization of meromorphic solutions of the equation $X^3 + Y^3 = 1$ we determine the zeros
of the function $\wp'(z) \pm \sqrt{3}$.
\end{abstract}

\section{Foreword}
Among the planar lattices two are singled out for their special symmetries: The square lattice and the hexagonal or triangular lattice.

In this short article we discuss certain properties of the Weierstrass $\wp$-function whose period lattice is hexagonal.

\section{Preliminaries}
The general facts about elliptic functions used in this section can be found in the books \cite{A}, \cite{C}, \cite{D}, \cite{H}, \cite{L}, and \cite{M-M}.

Let $\wp(z)$ be the Weierstrass $\wp$-function with $g_2 = 0$ and $g_3 = 1$, namely the $\wp$-function satisfying
\begin{equation}
\wp'(z)^2 = 4 \wp(z)^3 - 1
= 4 \left[\wp(z) - e_1\right]\left[\wp(z) - e_2\right]\left[\wp(z) - e_3\right]
\label{A5}
\end{equation}
(with $\lim_{z \to 0} z^2 \wp(z) = 1$),
where
\begin{equation}
e_1 = \frac{1}{4^{1/3}},
\qquad
e_2 = \frac{1}{4^{1/3}} \, e^{-2\pi i/3}
\qquad \text{and}\qquad
e_3 = \frac{1}{4^{1/3}} \,  e^{2\pi i/3}.
\label{A5a}
\end{equation}

If we set
\begin{equation}
\wp_k(z) := k^2 \wp(kz),
\qquad
k \in \mathbb{C} \smallsetminus \{0\},
\label{B1}
\end{equation}
then $\wp_k(z)$ is the Weierstrass $\wp$-function satisfying
\begin{equation}
\wp_k'(z)^2 = 4 \wp_k(z)^3 - k^6.
\label{B2}
\end{equation}
By choosing $k = e^{\pi i/3}$ we see that $\wp(z)$ possesses the symmetry
\begin{equation}
\wp\left(e^{\pi i/3} z\right) = e^{-2\pi i/3} \wp(z),
\qquad
\wp'\left(e^{\pi i/3} z\right) = -\wp'(z),
\label{B3}
\end{equation}
in particular, $\wp(z)$ is even and $\wp'(z)$ is odd.

Now, in view of \eqref{A5}, one period of $\wp(z)$ is
\begin{align}
\varpi &= \int_{4^{-1/3}}^{\infty} \frac{dx}{\sqrt{x^3 - 1/4}}
= \frac{2^{1/3}}{3} \int_0^1 \xi^{-5/6} (1 - \xi)^{-1/2} d\xi
\nonumber
\\
&= \frac{2^{1/3} \Gamma(1/6) \sqrt{\pi}}{3 \, \Gamma(2/3)}
= \frac{1}{2 \pi}\, \Gamma(1/3)^3
\label{B4}
\end{align}
(to obtain the last equation we used Legendre's duplication formula for the Gamma function \cite{A}, \cite{C}, together with the fact that
$\Gamma(1/3)\, \Gamma(2/3) = 2\pi/\sqrt{3}$). Then, with the help of \eqref{B3} we can see that
\begin{equation}
\text{the quantity} \quad e^{\pi i/3} \varpi \quad \text{is also a period of }\; \wp(z)
\label{B4a}
\end{equation}
and that $\varpi$ and $e^{\pi i/3} \varpi$ are, actually, primitive periods of $\wp(z)$. In fact, in view of \eqref{A5}, \eqref{A5a}, \eqref{B4}, and \eqref{B3}, we have
\begin{equation}
\wp\left(\frac{\varpi}{2}\right) = \frac{1}{4^{1/3}}
\quad\text{and}\quad
\wp\left(\frac{e^{\pi i/3}\varpi}{2}\right) = e^{-2\pi i/3} \wp\left(\frac{\varpi}{2}\right) = \frac{1}{4^{1/3}} \, e^{-2\pi i/3}.
\label{B4b}
\end{equation}

It follows that the period lattice of $\wp(z)$ is the hexagonal lattice
\begin{equation}
\mathbb{T} := \{m \varpi + n e^{\pi i/3} \varpi \, : \, m, n \in \mathbb{Z} \},
\label{B5}
\end{equation}
where
$\varpi$ is given by \eqref{B4}. In other words
\begin{equation}
\mathbb{T} = \varpi \mathbb{Z}[e^{\pi i/3}] = \varpi \mathbb{Z}[e^{2\pi i/3}]
\label{B5e}
\end{equation}
(the second equality follows from the fact that $e^{2\pi i/3} = e^{\pi i/3} - 1$), where $\mathbb{Z}[e^{2\pi i/3}]$ is the ring of Eisenstein
integers. Notice that $\mathbb{T}$ possesses the rotational symmetry
\begin{equation}
e^{\pi i/3} \mathbb{T} = \left\{e^{\pi i/3} \omega \,:\, \omega \in \mathbb{T}\right\} = \mathbb{T}
\label{B5s}
\end{equation}
and (since, e.g., $e^{-\pi i/3} \mathbb{T} = \mathbb{T}$) the symmetry with respect to complex conjugation, i.e.
\begin{equation}
\mathbb{T}^{\ast} := \left\{\bar{\omega} \,:\, \omega \in \mathbb{T}\right\}  = \mathbb{T}
\label{B5c}
\end{equation}
(as usual, $\bar{\omega}$ denotes the complex conjugate of $\omega$). Finally, the fundamental cell of $\mathbb{T}$ is
\begin{equation}
\mathcal{K} = \{s \varpi + t e^{\pi i/3} \varpi \, : \, 0 \leq s,t <1\}.
\label{B5a}
\end{equation}

From now on, unless otherwise stated, $\wp(z)$ will denote the specific Weierstrass $\wp$-function whose period lattice is $\mathbb{T}$.

We continue by considering the Weierstrass $\sigma$-function and the Weierstrass $\zeta$-function associated to $\mathbb{T}$, defined as
\begin{equation}
\sigma(z) = \sigma(z; \mathbb{T}) := z \prod_{\omega \in \mathbb{T}^{'}} \left(1 - \frac{z}{\omega}\right) e^{z/\omega + \frac{1}{2}(z/\omega)^2},
\qquad \text{where }\;
\mathbb{T}^{'} := \vec{\mathbb{T}} \smallsetminus \{0\},
\label{B6}
\end{equation}
and
\begin{equation}
\zeta(z) = \zeta(z; \mathbb{T}) := \frac{\sigma'(z)}{\sigma(z)}
= \frac{1}{z} + \sum_{\omega \in \mathbb{T}^{'}} \left(\frac{1}{z - \omega} + \frac{1}{\omega} + \frac{z}{\omega^2}\right).
\label{B6z}
\end{equation}
For typographical convenience the $\sigma$- and the $\zeta$-functions associated to the specific lattice $\mathbb{T}$ of \eqref{B5}
will be denoted by $\sigma(z)$ and $\zeta(z)$ respectively, instead of $\sigma(z; \mathbb{T})$ and $\zeta(z; \mathbb{T})$.

The function $\sigma(z)$ is entire of order $2$ and it is obvious from \eqref{B6} that its zeros are simple and located at the points of
$\mathbb{T}$. In particular, at $z=0$ we have
\begin{equation}
\sigma(0) = 0
\qquad\text{and}\qquad
\sigma'(0) = 1.
\label{B6b}
\end{equation}
Furthermore, from the symmetries \eqref{B5s} and \eqref{B5c} of $\mathbb{T}$ it follows that
\begin{equation}
\sigma\left(e^{\pi i/3} z\right) = e^{\pi i/3} \sigma(z)
\qquad\text{and}\qquad
\overline{\sigma(z)} = \sigma\left(\bar{z}\right),
\label{B6a}
\end{equation}
in particular, $\sigma(z)$ is odd and if $x$ denotes a real variable, then $\sigma(x)$ is real.

The function $\zeta(z)$ is meromorphic with simple poles located at the points of $\mathbb{T}$, and from \eqref{B6z} and \eqref{B6a} we obtain that
it satisfies the symmetry relations
\begin{equation}
\zeta\left(e^{\pi i/3} z\right) = e^{-\pi i/3} \zeta(z)
\qquad\text{and}\qquad
\overline{\zeta(z)} = \zeta\left(\bar{z}\right),
\label{B6s}
\end{equation}
in particular, $\zeta(z)$, too, is odd. Furthermore, in each of the intervals
\begin{equation}
I_n := (n\varpi, n\varpi + \varpi),
\qquad
n \in \mathbb{Z},
\label{BZ0}
\end{equation}
$\zeta(x)$ is decreasing with $\zeta(n\varpi+) = +\infty$ and $\zeta((n+1)\varpi-) = -\infty$.

Recall that, having $\sigma(z)$ and $\zeta(z)$, the function $\wp(z)$ can be constructed as
\begin{equation}
\wp(z) = -\zeta'(z) = \frac{\sigma'(z)^2 - \sigma''(z) \sigma(z)}{\sigma(z)^2}
= \frac{1}{z^2} + \sum_{\omega \in \mathbb{T}^{'}} \left[\frac{1}{(z - \omega)^2} - \frac{1}{\omega^2}\right]
\label{B7c}
\end{equation}
and from \eqref{B7c} it follows immediately that
\begin{equation}
\wp'(z) = -\sum_{\omega \in \mathbb{T}} \frac{2}{(z - \omega)^3}.
\label{B7d}
\end{equation}
We remind the reader that the differential equation \eqref{A5} is obtained by noticing that its left- and right-hand sides are elliptic functions with matched poles. As for the symmetry relations \eqref{B3}, let us notice that they also follow directly from the first equality of \eqref{B6s}.

The function $\wp'(z)$ has three zeros in $\mathcal{K}$ and the fact that it is odd and $\mathbb{T}$-periodic implies that these zeros are
$\varpi/2$, $e^{\pi i/3} \varpi/2$, and $(1 + e^{\pi i/3}) \varpi/2 = \sqrt{3}\, e^{\pi i/6} \varpi/2$.

If we diferentiate both sides of \eqref{A5} and eliminate $\wp'(z)$, we obtain the equation
\begin{equation}
\wp''(z) = 6 \wp(z)^2,
\label{B7dd}
\end{equation}
thus, in the interval, say, $(0, \varpi)$ the function $\wp(x)$ is convex, with $\wp(0+) = \wp(\varpi-) = +\infty$,
while $\wp'(x)$ is increasing, with $\wp'(0+) = -\infty$ and
$\wp'(\varpi-) = +\infty$. And since $\wp'(\varpi/2) = 0$, the minimum value of $\wp(x)$ is $\wp(\varpi/2) = 4^{-1/3}$.

In connection to the periods $\varpi$ and $e^{\pi i/3} \varpi$ of $\wp(z)$, it is expedient to introduce the quantities
\begin{equation}
\eta_1 := \zeta(z + \varpi) - \zeta(z)
\qquad\text{and}\qquad
\eta_2 := \zeta\left(z + e^{\pi i/3} \varpi\right) - \zeta(z).
\label{BZ1}
\end{equation}
From the fact that the derivative $\zeta'(z)$ is elliptic with periods $\varpi$ and $e^{\pi i/3} \varpi$ we infer that $\eta_1$ and
$\eta_2$ are independent of $z$. By setting $z = -\varpi/2$ in the first equation and $z = -e^{\pi i/3}\varpi/2$ in the second equation of \eqref{BZ1} we obtain, in view of \eqref{B6s},
\begin{equation}
\eta_1 = 2 \,\zeta\left(\frac{\varpi}{2}\right)
\qquad\text{and}\qquad
\eta_2 = 2 \,\zeta\left(\frac{e^{\pi i/3}\varpi}{2}\right) = 2 e^{-\pi i/3} \zeta\left(\frac{\varpi}{2}\right)
= e^{-\pi i/3} \eta_1.
\label{BZ2}
\end{equation}
Now, a well-known property of $\wp(z)$ is that $\eta_1$ and $\eta_2$ satisfy \textit{Legendre's relation}
\begin{equation}
\eta_1 e^{\pi i/3}\varpi - \eta_2 \varpi = 2 \pi i
\label{BZ3}
\end{equation}
(which follows by integrating $\zeta(z)$ along the perimeter of the parallelogram with vertices $a$, $a + \varpi$, $a + e^{\pi i/3}\varpi$, and
$a + \varpi + e^{\pi i/3}\varpi$, where $a$ is any complex number in the interior of $\mathcal{K}$).

By combining \eqref{BZ2} and \eqref{BZ3} we obtain
\begin{equation}
2 \varpi \zeta\left(\frac{\varpi}{2}\right) = \eta_1 \varpi = \frac{2 \pi}{\sqrt{3}},
\label{BZ4}
\end{equation}
thus, in view of \eqref{B4},
\begin{equation}
\eta := \eta_1 = 2 \,\zeta\left(\frac{\varpi}{2}\right) = \frac{2 \pi}{\sqrt{3}}\, \frac{1}{\varpi}
= \frac{4 \pi^2}{\sqrt{3}} \, \frac{1}{\Gamma(1/3)^3}
\label{BZ5}
\end{equation}
(from now on, for typographical convenience we will write $\eta$ instead of $\eta_1$).

Incidentally, by using \eqref{BZ4} in \eqref{B6z} we get a little bonus:
\begin{equation}
\frac{2\pi}{\sqrt{3}} = 2\varpi \zeta\left(\frac{\varpi}{2}\right)
= 4 + \sum_{\kappa \in \mathbb{Z}[e^{\pi i/3}] \smallsetminus \{0\}} \frac{1}{(1- 2\kappa) \kappa^2}
\label{BZ6}
\end{equation}
(if instead of the hexagonal lattice and the Eisenstein integers we consider the square lattice and its corresponding Weierstrass $\zeta$-function,
then, by following the same approach, we can obtain the value of the above sum taken over the nonzero Gaussian integers).

Finally, in view of \eqref{B6z}, the first equation in \eqref{BZ1} can be written as
\begin{equation*}
\eta = \frac{\sigma'(z + \varpi)}{\sigma(z + \varpi)} - \frac{\sigma'(z)}{\sigma(z)},
\end{equation*}
which implies
\begin{equation*}
\sigma(z + \varpi) = C e^{\eta z} \sigma(z),
\end{equation*}
where $C$ is a constant. Setting $z = - \varpi/2$ and recalling that $\sigma(z)$ is odd we get $C = -e^{\eta \varpi/2}$ or, in view of
\eqref{BZ4}, $C = -e^{\pi/\sqrt{3}}$, thus
\begin{equation}
\sigma(z + \varpi) = -e^{\pi/\sqrt{3}} e^{\eta z} \sigma(z)
\label{BZ7}
\end{equation}
(since $\varpi, \eta >0$, formula \eqref{BZ7} implies that $\lim_{x \to -\infty} \sigma(x) = 0$, while $\limsup_{x \to +\infty} \sigma(x) = +\infty$
and $\liminf_{x \to +\infty} \sigma(x) = -\infty$).

In the same manner, the second equation in \eqref{BZ1} yields
\begin{equation}
\sigma\left(z + e^{\pi i/3}\varpi\right) = -e^{\pi/\sqrt{3}} \exp\left(e^{-\pi i/3}\eta z\right) \sigma(z).
\label{BZ8}
\end{equation}

\section{The zeros of $\wp(z)$}
The following sets will be useful in the sequel:
\begin{equation}
\mathcal{T}_1 := \text{the equilateral triangle with vertices }\; 0, \ \varpi, \ e^{\pi i/3} \varpi;
\label{T1}
\end{equation}
\begin{equation}
\mathcal{T}_2 := \text{the equilateral triangle with vertices }\; \varpi, \ e^{\pi i/3} \varpi, \
\sqrt{3}\, e^{\pi i/6} \varpi
\label{T2}
\end{equation}
(notice that $\sqrt{3}\, e^{\pi i/6} \varpi = \varpi + e^{\pi i/3} \varpi$);
\begin{equation}
\mathcal{K}\,' := \{s \varpi + t e^{2\pi i/3} \varpi \, : \, 0 \leq s,t <1\}
\label{C17}
\end{equation}
the cell formed by the primitive periods $\varpi$ and $e^{2\pi i/3} \varpi$.
Notice that the union of $\mathcal{T}_1$ and $\mathcal{T}_2$ is the fundamental cell $\mathcal{K}$ of \eqref{B5a}, while the intersection of
$\mathcal{K}\,'$ with $\mathcal{K}$ is the equilateral triangle $\mathcal{T}_1$.

The following proposition is a special case of a result regarding the location of the zeros of the general Weierstrass $\wp$-function, which can be
found in \cite{D}. We include a proof here for the sake of completeness.

\medskip

\textbf{Proposition 1.} The zeros of $\wp(z)$ are
\begin{equation}
\pm r \varpi + \omega,
\qquad
\omega \in \mathbb{T},
\qquad \text{where }\;
r := \frac{\sqrt{3}}{3} \, e^{\pi i/6}.
\label{C1b}
\end{equation}
Furthermore, they are all simple.

\smallskip

\textit{Proof}. The zeros of $\wp(z)$ are simple, since if $\wp(z^{\star}) = 0$, then \eqref{A5} implies that $\wp'(z^{\star}) \ne 0$. We, also, know that $\wp(z)$ has exactly two zeros in each period cell. Furthermore, as we have seen, $\wp(x) > 0$ for $x \in (0, \varpi)$. Therefore, by
\eqref{B3} we get that $\wp(x e^{\pi i/3}) \ne 0$ and consequently $\wp(z) \ne 0$ for $z \in \partial \mathcal{K} \cup \partial\mathcal{K}\,'$. Thus,
the two zeros of $\wp(z)$ in $\mathcal{K}$ lie in the union of the interiors of the two equilateral triangles $\mathcal{T}_1$
and $\mathcal{T}_2$ of \eqref{T1} and \eqref{T2}.

Suppose $z^{\star}$ is a zero of $\wp(z)$ in the interior of $\mathcal{T}_1$. Then, since $\wp(z)$ is even and $\mathbb{T}$-periodic,
$0 = \wp(-z^{\star}) = \wp(\varpi + e^{\pi i/3} \varpi -z^{\star})$. Thus, $\varpi + e^{\pi i/3} \varpi -z^{\star}$ is a zero of $\wp(z)$
lying in the interior of $\mathcal{T}_2$  (likewise, if $z^{\star}$ is a zero in the interior of $\mathcal{T}_2$, then
$\varpi + e^{\pi i/3} \varpi -z^{\star}$ is a zero of in the interior of $\mathcal{T}_1$). It follows that $\wp(z)$ has one zero, say $z_1$
in the interior of $\mathcal{T}_1$ and one zero, say $z_2$ in the interior of $\mathcal{T}_2$, and these are the only zeros of $\wp(z)$ in
$\mathcal{K}$.

Now by \eqref{B3} we get that $e^{\ell\pi i/3} z_1$, $\ell = 0, 1, \ldots, 5$, are zeros of $\wp(z)$. These zeros lie at the vertices of a regular hexagon centered at $0$. Hence, the $\mathbb{T}$-periodicity of $\wp(z)$ implies that for every point $\omega \in \mathbb{T}$
we have an associated set of zeros $\mathcal{Z}(\omega) := \{\omega + e^{\ell\pi i/3} z_1 \,:\, \ell = 0, 1, \ldots, 5\}$ lying at the vertices of
a regular hexagon centered at $\omega$.

From the above observations it follows that
\begin{equation}
\mathcal{Z}(0) \cap \mathcal{T}_1 = \mathcal{Z}(\omega) \cap \mathcal{T}_1
= \mathcal{Z}(e^{\pi i/3}\omega) \cap \mathcal{T}_1 = \{z_1\}
\label{E1}
\end{equation}
and, consequently, $z_1$ is equidistant from the vertices $0$, $\omega$, and $e^{\pi i/3}\omega$ of $\mathcal{T}_1$. Therefore, $z_1$ is located at
the center of $\mathcal{T}_1$ and, likewise $z_2$ is located at the center of $\mathcal{T}_2$. In other words, the zeros of $\wp(z)$ in $\mathcal{K}$ are
\begin{equation}
\frac{\sqrt{3}}{3} \, e^{\pi i/6} \varpi
\qquad\text{and}\qquad
\frac{2 \sqrt{3}}{3} \, e^{\pi i/6} \varpi
\label{C1a}
\end{equation}
and the proof is completed by noticing that
\begin{equation}
\frac{2 \sqrt{3}}{3} \, e^{\pi i/6} \varpi = -\frac{\sqrt{3}}{3} \, e^{\pi i/6} \varpi
\quad(\text{mod } \mathbb{T})
\label{C2}
\end{equation}
(due to the fact that $\sqrt{3}\, e^{\pi i/6} \varpi = (1 + e^{\pi i/3}) \varpi \in \mathbb{T}$).
\hfill $\blacksquare$

\medskip

The knowledge of its zeros and its poles yields an expression of $\wp(z)$ in terms of the $\sigma$-function, which is quite different from the one of
\eqref{B7c} \cite{A}, \cite{L}
\begin{equation}
\wp(z) = -\frac{1}{\sigma\left(r\varpi\right)^2} \,
\frac{\sigma\left(z - r\varpi\right) \sigma\left(z + r\varpi\right)}{\sigma(z)^2}.
\label{C3}
\end{equation}
In the special case studied in this article, where the period lattice is the hexagonal lattice $\mathbb{T}$, there is, yet, one more expression of
$\wp(z)$ in terms of $\sigma(z)$.

\medskip

\textbf{Proposition 2.} Let $\wp(z)$ be the $\mathbb{T}$-periodic Weierstrass $\wp$-function, where $\mathbb{T}$ is the hexagonal lattice of
\eqref{B5}, and $\sigma(z)$ the associated $\sigma$-function. Then,
\begin{equation}
\wp(z) = \frac{\sigma(z; r \mathbb{T})}{\sigma(z)^3}
= \frac{r \sigma\left(r^{-1} z\right)}{\sigma(z)^3},
\label{C9}
\end{equation}
where $r$ is given by \eqref{C1b}.

\smallskip

\textit{Proof}. First, let us observe that
\begin{equation}
\mathbb{T} \, \cup \left(\mathbb{T} + r\varpi\right)
\cup \left(\mathbb{T} - r\varpi\right)
= r \mathbb{T},
\label{C4}
\end{equation}
from which it follows that the functions
\begin{equation}
\sigma(z) \sigma(z - r\varpi) \sigma(z + r\varpi)
\quad\text{and}\quad
\sigma(z; r \mathbb{T}) = r \sigma\left(r^{-1} z\right)
\label{C5}
\end{equation}
have the same zeros (including multiplicities).
Both functions in \eqref{C5} are odd and entire of order $2$, while their derivatives at $z = 0$ are $-\sigma(r\varpi)^2$ and $1$ respectively.
Therefore, there is a constant $a \in \mathbb{C}$ such that
\begin{equation}
\sigma(z) \sigma(z - r\varpi) \sigma(z + r\varpi)
= - r \sigma(r\varpi)^2 e^{a z^2} \sigma\left(r^{-1} z\right)
\label{C8}
\end{equation}
and, consequently, \eqref{C3} implies
\begin{equation}
\wp(z) = \frac{r \sigma\left(r^{-1} z\right)}{\sigma(z)^3} \, e^{a z^2}
\label{C10}
\end{equation}
(actually, the symmetries $\bar{r} \mathbb{T} = r e^{-\pi i/3} \mathbb{T} = r \mathbb{T} = r \mathbb{T}^{\ast}$ imply that
$a \in \mathbb{R}$).
Since $\wp(z + \varpi) = \wp(z)$, formula \eqref{C10} yields
\begin{equation}
\frac{\sigma\left(r^{-1} z + r^{-1} \varpi\right)}{\sigma(z + \varpi)^3} \, e^{2az\varpi + a\varpi^2}
= \frac{\sigma\left(r^{-1} z\right)}{\sigma(z)^3}.
\label{C11}
\end{equation}
Noticing that
\begin{equation}
r^{-1} = \sqrt{3} \, e^{\pi i/6}  = e^{-\pi i/3} + 1
\label{C12}
\end{equation}
we have, in view of \eqref{BZ7} and \eqref{BZ4},
\begin{align}
\sigma\left(r^{-1} z + r^{-1} \varpi\right) &= \sigma\left(r^{-1} z + e^{-\pi i/3} \varpi  + \varpi\right)
\nonumber
\\
&= - e^{\pi/\sqrt{3}} \exp\left(r^{-1} \eta z + \frac{2\pi}{\sqrt{3}}\, e^{-\pi i/3}\right)\sigma\left(r^{-1} z + e^{-\pi i/3} \varpi\right).
\label{C13}
\end{align}
Now, in view of \eqref{B6a} and \eqref{BZ8},
\begin{align}
\sigma\left(r^{-1} z + e^{-\pi i/3} \varpi\right) &= \overline{\sigma\left(\overline{r^{-1} z} + e^{\pi i/3} \varpi\right)}
\nonumber
\\
&= - e^{\pi/\sqrt{3}}\, \overline{\exp\left(e^{-\pi i/3} \, \eta \overline{r^{-1} z}\right)} \,
\overline{\sigma\left(\overline{r^{-1} z}\right)}
\nonumber
\\
&= - e^{\pi/\sqrt{3}}\, \exp\left(e^{\pi i/3} \, r^{-1} \eta z\right) \, \sigma\left(r^{-1} z\right).
\label{C14}
\end{align}
Substituting \eqref{C14} in \eqref{C13} gives
\begin{align}
\sigma\left(r^{-1} z + r^{-1} \varpi\right)
&= e^{2\pi/\sqrt{3}} \exp\left(\left(1 + e^{\pi i/3}\right) r^{-1} \eta z + \frac{2\pi}{\sqrt{3}}\, e^{-\pi i/3}\right) \,
\sigma\left(r^{-1} z\right)
\nonumber
\\
&= e^{2\pi/\sqrt{3}} \exp\left(3 \eta z + \frac{\pi}{\sqrt{3}} + i\pi\right) \,
\sigma\left(r^{-1} z\right)
\nonumber
\\
&= -e^{3\pi/\sqrt{3}} e^{3 \eta z} \sigma\left(r^{-1} z\right).
\label{C15}
\end{align}
Also, again by \eqref{BZ7}
\begin{equation}
\sigma(z + \varpi)^3 = -e^{3\pi/\sqrt{3}} e^{3 \eta z} \sigma(z)^3
\label{C16}
\end{equation}
and, finally, by substituting \eqref{C15} and \eqref{C16} in \eqref{C11} we get that $a = 0$.
\hfill $\blacksquare$

\medskip

Formula \eqref{C9} is reminiscent of the well-known equation \cite{C}
\begin{equation}
\wp'(z) = -\frac{\sigma(2z)}{\sigma(z)^4},
\label{B9a}
\end{equation}
which is valid for any lattice, not just $\mathbb{T}$.

\section{The zeros of $\wp'(z) \pm \sqrt{3}$}
Regarding the Weierstrass $\wp$-function associated to any period lattice, a consequence of the addition theorem is that if $n$ is an integer, then $\wp(nz)$ can be expressed as a rational function of $\wp(z)$. Thus, for any rational number $q$, the quantity $\wp(qz)$ is an algebraic function of
$\wp(z)$ and the same is true for the derivatives $\wp^{(k)}(qz)$. For instance,
\begin{align}
\wp\left(\frac{z}{2}\right) = \wp(z) &+ \sqrt{\left[\wp(z) - e_1\right] \left[\wp(z) - e_2\right]}
\nonumber
\\
&+ \sqrt{\left[\wp(z) - e_2\right] \left[\wp(z) - e_3\right]} + \sqrt{\left[\wp(z) - e_3\right] \left[\wp(z) - e_1\right]}
\label{G1}
\end{align}

In particular, if $\omega$ is a period of $\wp(z)$ and $g_2$, $g_3$ are algebraic numbers, then for any $q \in \mathbb{Q}$ and any
$k = 0, 1, 2, \ldots$ the number $\wp^{(k)}(qz)$ is algebraic or $\infty$.

In this section we will determine the values of $\wp(z)$ and $\wp'(z)$ for certain arguments $z$ without following the above general approach.

We start with a little detour. It is well known that the algebraic curve
\begin{equation}
X^3 + Y^3 = 1
\label{A4}
\end{equation}
can be ``uniformized" by elliptic functions (this is a consequence of the fact that the genus of the curve is $1$).
In fact, by using \eqref{A5} and the evenness of $\wp(z)$ it is easy to check that
\begin{equation}
X = f(z) := \frac{\wp'(z) + \sqrt{3}}{2\sqrt{3} \, \wp(z)}
\qquad
\text{and}
\qquad
Y = f(-z)
\label{A6}
\end{equation}
is a uniformization (i.e. global parametrization) of \eqref{A4}.

Incidentally, the function $f(z)$ above is related to the Dixon elliptic functions \cite{Di} which, too, unifomize the curve \eqref{A4}.

Formula \eqref{A6} can be viewed as a meromorphic (in $\mathbb{C}$) solution $(X, Y)$ of the Diophantine-type equation \eqref{A4}. And, as it turns
out, this solution is unique in the following sense.

\medskip

\textbf{Theorem} (I.N. Baker \cite{B}, 1966)\textbf{.} Every meromorphic solution of \eqref{A4} is of the form
\begin{equation}
X = f\big(h(z)\big),
\qquad
Y = \rho f\big(-h(z)\big),
\qquad\quad
\rho^3 = 1,
\label{A7}
\end{equation}
where $f(z)$ is as in \eqref{A6} and $h(z)$ is an entire function.

\medskip

Here is a brief sketch of Baker's neat and instructive argument: The function
$f(z)$ of \eqref{A6} is elliptic and has exactly three simple poles in each period cell; hence it takes every value there three times, counting
multiplicities. Since $[1 - f(z)^3]^{1/3} = \rho f(-z)$ is meromorphic it follows that $f(z) - e^{2\pi k i/3}$ has a triple zero for
$k = 0, 1, 2$. Thus, in any fixed period cell, for each such $k$ there is a unique $z_k$  for which $f(z_k) = e^{2\pi k i/3}$. Now $f'(z)$ has
exactly three double poles in that period cell, hence it takes every value there exactly six times. In particular, $z_k$, $k = 0, 1, 2$, are double
zeros of $f'(z)$ and, consequently, these are the only zeros of $f'(z)$ in that period cell. Since the singularities of the (multi-valued) inverse
function $f^{-1}(w)$ arise at the values $w = f(z)$ for which $f'(z) = 0$, it follows that the singularities of $f^{-1}(w)$ lie over
$w_s = 1, e^{2\pi i/3}, e^{-2\pi i/3}$. From the above observations it is not hard to show that if $X(z)$, $Y(z)$ is any meromorphic solution of \eqref{A4}, then the function element $h(z) = f^{-1}(X(z))$ can be analytically continued indefinitely along any curve $\gamma$ in the complex plane,
even if $\gamma$ passes through a point $z_s$ for which $X(z_s) = w_s$, since the fact that $[1 - X(z)^3]^{1/3}$ is meromorphic implies that
$X(z) = w_s + \phi(z)^3$, where $\phi(z)$ is analytic near $z_s$ and $\phi(z_s) = 0$. Therefore, by the monodromy theorem $h(z)$ is an entire function.

\medskip

Baker's theorem has a remarkable implication regarding the values of $\wp'(z)$.

\medskip

%

\textbf{Proposition 3.} The zeros of $\wp'(z)+ \sqrt{3}$ are
\begin{equation}
\frac{\varpi}{3} + \omega,
\quad
e^{\pi i/3} \frac{2\varpi}{3} + \omega,
\quad
e^{2\pi i/3} \frac{\varpi}{3} + \omega,
\qquad
\omega \in \mathbb{T},
\label{D1}
\end{equation}
(clearly, they are all simple).

\smallskip

\textit{Proof}.
As we have seen, $\wp'(x)$ is increasing in
$(0, \varpi)$, with $\wp'(0+) = -\infty$ and
$\wp'(\varpi/2) = 0$. Thus, there is a unique $x_+ \in (0, \varpi)$ such that $\wp'(x_+) = -\sqrt{3}$, and we must have $x_+ < \varpi/2$.
From the fact that $\wp'(z)$ is odd and $\varpi$-peridic it follows that if $x_- := \varpi - x_+$, then $\wp'(x_-) = \sqrt{3}$.

Each of the functions $\wp'(z) + \sqrt{3}$ and $-\wp'(z) + \sqrt{3}$ has exactly three zeros in $\mathcal{K}\,'$ (counting multiplicities).
Therefore, the second equation in \eqref{B3} implies that the zeros of $\wp'(z) + \sqrt{3}$ are
\begin{equation}
\omega + x_+,
\quad
\omega + e^{\pi i/3} x_-,
\quad
\omega + e^{2\pi i/3} x_+,
\qquad
\omega \in \mathbb{T},
\label{C18a}
\end{equation}
and it remains to show that $x_+ = \varpi/3$.

From the fact that $X$ and $Y$ of \eqref{A6} satisfy \eqref{A4} we obtain that
\begin{equation*}
\left[\wp'(z) + \sqrt{3}\right]^3 + \left[\wp'(-z) + \sqrt{3}\right]^3 = \left[2\sqrt{3} \, \wp(z)\right]^3,
\end{equation*}
which, by dividing by $\left[\wp'(z) + \sqrt{3}\,\right]^3$ implies (since $\wp'(z)$ is odd)
\begin{equation}
\left[\frac{2\sqrt{3} \, \wp(z)}{\wp'(z) + \sqrt{3}}\right]^3 + \left[\frac{\wp'(z) - \sqrt{3}}{\wp'(z) + \sqrt{3}}\right]^3 = 1.
\label{F1}
\end{equation}
Since \eqref{F1} yields another solution of \eqref{A4}, Baker's theorem implies that there is an entire function $h(z)$ such that
\begin{equation}
\rho \frac{2\sqrt{3} \, \wp(z)}{\wp'(z) + \sqrt{3}} =
\frac{\wp'\big(h(z)\big) + \sqrt{3}}{2\sqrt{3} \, \wp\big(h(z)\big)}.
\label{F2}
\end{equation}
The quantity in the left-hand side of \eqref{F2} is a $\mathbb{T}$-periodic elliptic function with exactly three simple poles in $\mathcal{K}$.
Therefore, its period lattice must be $\mathbb{T}$ (i.e. it cannot be larger than $\mathbb{T}$). Thus, the same must be true for the quantity in the
right-hand side of \eqref{F2}. Consequently,
\begin{equation}
h(z) = \alpha z + b,
\qquad\text{where }\; |\alpha| = 1.
\label{F3}
\end{equation}
Furthermore, \eqref{F2} also implies that the functions $\wp'\big(h(z)\big) + \sqrt{3}$ and $\sigma(z) \wp'(z)$ have the same zeros (including
multiplicities). Thus, in view of \eqref{C9}, the set of zeros of $\wp'\big(h(z)\big) + \sqrt{3}$ is the hexagonal lattice
\begin{equation}
\frac{e^{\pi i/6}}{\sqrt{3}} \, \mathbb{T}
\label{F4}
\end{equation}
(and, of course, all these zeros are simple).

Since the map $z \mapsto \alpha z + b$ preserves (Euclidean) distances in the complex plane, the matching of the zeros of
$\wp'\big(h(z)\big) + \sqrt{3}$ and $\sigma(z) \wp'(z)$ implies that the distance between the zeros $x_+$ and $e^{\pi i/3} x_-$ of
$\wp'\big(h(z)\big) + \sqrt{3}$ must be equal to the distance of neighboring zeros of $\sigma(z) \wp'(z)$, namely $\varpi/\sqrt{3}$
(in view of \eqref{F4}). Thus, since $x_- = \varpi - x_+$,
\begin{equation*}
\left|x_+ - e^{\pi i/3} (\varpi - x_+)\right| = \frac{\varpi}{\sqrt{3}},
\qquad\text{i.e.}\qquad
\frac{(3 x_+ - \varpi)^2}{4} + \frac{3(\varpi - x_+)^2}{4}= \frac{\varpi^2}{3},
\end{equation*}
which yields $x_+ = \varpi/3$ or $x_+ = 2\varpi/3$. Therefore, $x_+ = \varpi/3$ since, as we have seen, $x_+ < \varpi/2$.
\hfill $\blacksquare$

\medskip

An immediate consequence of Proposition 3 is that the zeros of $\wp'(z) - \sqrt{3}$ are
\begin{equation}
\frac{2\varpi}{3} + \omega,
\quad
e^{\pi i/3} \frac{\varpi}{3} + \omega,
\quad
e^{2\pi i/3} \frac{2\varpi}{3} + \omega,
\qquad
\omega \in \mathbb{T},
\label{Fe}
\end{equation}
(and, of course, they are all simple).

Finally, let us notice that the equation \eqref{A5} implies
\begin{equation}
\wp'(z)^2 - 3 = 4 \left[\wp(z)^3 - 1\right]
\label{C20}
\end{equation}
from which we get that (since $\wp(x)$ is real for $0 < x < \varpi$)
\begin{equation}
\wp\left(\pm\frac{\varpi}{3}\right) = 1
\label{C21}
\end{equation}
or, equivalently,
\begin{equation}
\int_1^{\infty} \frac{dx}{\sqrt{4x^3 - 1}} = \int_{\wp\left(\frac{\varpi}{3}\right)}^{\infty} \frac{dx}{\sqrt{4x^3 - 1}}
= \frac{\varpi}{3}
= \frac{1}{6 \pi}\, \Gamma(1/3)^3,
\label{C22}
\end{equation}
where the last equality is obtained from \eqref{B4}.

%
%
%
%
%
%
%
%
%
%
%
%

\end{document}